\documentclass[10pt,twoside]{classe-ma_e-ppt}
\usepackage{amssymb,amsbsy,amsmath,amsfonts,amssymb}
\usepackage[english,francais]{babel}
\usepackage{times}
\newtheorem{theorem}{Theorem}[section]
\newtheorem{theoreme}{Th\'eor\`eme}

\renewcommand{\Sect}[1]{\renewcommand{\thesection}{\Roman{section}}\section{#1}
\renewcommand{\thesection}{\arabic{section}.}
}

\newcommand{\R}{{\mathbb R}}
\newcommand{\N}{{\mathbb N}}
\newcommand{\fer}[1]{(\ref{#1})}
\newcommand{\be}{\begin{equation}}
\newcommand{\ee}{\end{equation}}
\newcommand{\dis}{\displaystyle}
\def\<#1>{\left\langle #1\right\rangle}
\def\lambdaper{\lambda_{\text{\rm per}}}
 
\AuteurCourant{J. Clairambault, S. Gaubert, B. Perthame}
\TitreCourant{Perron and Floquet eigenvalues}
\Recu{}{}{}
\begin{document}\selectlanguage{english}
\title{An inequality for the Perron and Floquet eigenvalues of monotone differential systems and age
structured equations}
\author{Jean Clairambault $^{\text{a,d}}$,\ \ St\'ephane Gaubert
$^{\text{b}}$, \ \  Beno\^\i t Perthame~$^{\text{a, c}}$}
\address{\begin{itemize}\labelsep=2mm\leftskip=-5mm
\item[$^{\text{a}}$]
INRIA, projet BANG,
Domaine de Voluceau, BP 105, F78153 Le Chesnay Cedex,\\ E-mail:
jean.clairambault@inria.fr \\
\item[$^{\text{b}}$]
INRIA, projet MAXPLUS,
Domaine de Voluceau, BP 105, F78153 Le Chesnay Cedex,\\ E-mail:
stephane.gaubert@inria.fr \\
\item[$^{\text{c}}$]
D\'epartement de Math\'ematiques et Applications, UMR 8553, \\
\'Ecole Normale Sup\'erieure, 45, rue d'Ulm, F75230 Paris Cedex 05,
\\Benoit.Perthame@ens.fr\\ 
\item[$^{\text{d}}$]
INSERM U 776 ``Rythmes biologiques et cancers'', H\^opital
Paul-Brousse,\\ 14, Av. Paul-Vaillant-Couturier, F94807 Villejuif
Cedex
\end{itemize}}
\renewcommand{\thefootnote}{}
\footnotetext{The second author was partially supported by the joint RFBR-CNRS grant number 05-01-02807}
\maketitle\thispagestyle{empty}
\begin{Abstract}{
For monotone linear differential systems with periodic coefficients, the
(first) Floquet eigenvalue measures the growth rate of the system. We define
an appropriate arithmetico-geometric time average of the coefficients
for which  we can prove that the  Perron eigenvalue is smaller
than the Floquet eigenvalue. We apply this method to Partial
Differential Equations, and we use it for an age-structured systems 
of equations for the cell cycle. This opposition between Floquet and 
Perron eigenvalues models the loss of circadian rhythms by cancer cells. 
}\end{Abstract}
\selectlanguage{french}
\begin{Ftitle}{Une in\'egalit\'e pour les valeurs propres de Floquet et de Perron de syst\`emes diff\'erentiels monotones et d'\'equations structur\'ees en \^age
} 
\end{Ftitle}
\begin{Resume}{
La (premi\`ere) valeur propre de Floquet d\'ecrit
le taux de croissance des syst\`emes diff\'erentiels
lin\'eaires monotones \`a coefficients p\'eriodiques. Nous d\'efinissons une 
moyenne arithm\'etico-g\'eom\'etrique
en temps des coefficients, qui nous permet de d\'emontrer que la valeur
propre de Perron pour le syst\`eme ainsi moyenn\'e est plus petite que celle de
Floquet. La m\'ethode s'applique aux \'Equations aux D\'eriv\'ees Partielles et nous l'utilisons
pour un syst\`eme d'\'equa\-tions structur\'ees en \^age qui d\'ecrit 
le cycle cellulaire.  Cette opposition entre valeurs propres de Floquet et de Perron mod\'elise la
perte de contr\^ole circadien pour le cycle cellulaire des cellules canc\'ereuses. 
}\end{Resume}

%
%
%

%
\AFv

Les syst\`emes biologiques sont souvent soumis \`a des contr\^oles
p\'eriodiques. Un exemple en est fourni par le
rythme circadien (journalier) qui trouve son origine au niveau des noyaux
suprachiasmatiques de l'hypothalamus dans un r\'eseau de r\'egulation 
g\'enique \`a pr\'esent bien \'etudi\'e
\cite{nagoshi, BernardPlos}. Une question m\'edicale reli\'ee est de
comprendre l'influence de ce rythme sur la croissance des populations de cellules tumorales, en
particulier dans une perspective th\'erapeutique. Un cadre
math\'ematique naturel pour \'etudier cette question est fourni par
la th\'eorie des \'equations physiologiquement structur\'ees. 

Dans cette note, nous proposons de comparer la (premi\`ere) valeur propre de
Floquet, qui d\'ecrit la croissance d'un syst\`eme dont les 
param\`etres sont soumis \`a un contr\^ole p\'eriodique, et la valeur 
propre de Perron, qui d\'ecrit la croissance du m\^eme syst\`eme, mais
\`a coefficients moyenn\'es. Pour parvenir \`a cette
comparaison, nous introduisons une moyenne
arithm\'etico-g\'eom\'etrique appropri\'ee, dans trois cas : 
\'Equations Diff\'erentielles Ordinaires, syst\`eme 
dynamique discret,  \'Equations aux
D\'eriv\'ees Partielles (EDP) structur\'ees en \^age.

Le premier cas concerne un syst\`eme diff\'erentiel lin\'eaire
monotone p\'eriodique. Soit $t\mapsto A(t)$ une application
$T$-p\'eriodique, \`a valeurs dans $\R^{d\times d}$,
int\'egrable sur $[0,T]$, et consid\'erons
l'\'equation diff\'erentielle $\dot X (t)= A(t)X(t)$, o\`u
$X$ est absolument continue et $X(0)$ est prescrit.
Nous supposerons que pour $i\neq j$, on a $A_{ij}(t)\geq 0$ pour presque
tout $t$, ce qui garantit que le flot en temps positif
associ\'e \`a cette \'equation diff\'erentielle pr\'eserve
l'ordre partiel usuel de $\R^d$. Les propri\'et\'es
spectrales de cette \'equation diff\'erentielle rel\`event alors
de la th\'eorie de Perron-Frobenius. En particulier, la premi\`ere
valeur propre de Floquet, $\lambdaper$, est le plus grand r\'eel
tel qu'il existe une fonction absolument continue
$X$, $T$-p\'eriodique, \`a valeurs
dans $\R_+^d$ (o\`u $\R_+$ d\'esigne l'ensemble
des r\'eels positifs ou nuls), et non identiquement nulle, telle que
l'\'equation diff\'erentielle~\eqref{eq:floquet} soit satisfaite.
Rappelons que la moyenne arithm\'etique d'une fonction $T$-p\'eriodique
$u(t)$ est donn\'ee par~\eqref{e-ar}, et  d\'efinissons
la matrice $\bar A$ par~\eqref{eq:agaver}.
Cette matrice \`a coefficients hors diagonaux positifs ou nuls
est dot\'ee d'une valeur propre de Perron classique $\lambda_s$,
qui est le plus grand r\'eel tel qu'il existe un vecteur non-nul
$U\in \R_+^d$ tel que $ \lambda_s U = \bar A U$.
Le r\'esultat suivant peut \^etre vu comme une 
g\'en\'eralisation de celui de \cite{clmp}.
\begin{theoreme}
On a toujours 
$\lambdaper  \geq \lambda_s $.
\end{theoreme}

Une in\'egalit\'e analogue est aussi obtenue pour les deux autres cas 
consid\'er\'es, d'un syst\`eme discret en temps et du syst\`eme d'EDP structur\'ees
en \^age  \fer{eq:cellcycle} qui  d\'ecrit le cycle de division cellulaire.

\par\medskip\centerline{\rule{2cm}{0.2mm}}\medskip
\setcounter{section}{0}\setcounter{equation}{0}
\selectlanguage{english}
\Sect{Introduction}

Biological systems are often subject to periodic controls. This
occurs for instance with circadian rhythms, the origin of which are 
found in the suprachiasmatic nuclei of the hypothalamus, in a now well established
gene regulatory network  \cite{nagoshi, BernardPlos}. A related medical
question is to understand the interactions between the cell cycle and
this circadian rhythm, which is expressed in every nucleated cell, 
with coordinating inputs from the hypothalamus. How can these rhythms induce differentiated
growth between healthy and tumoral cells? A molecular 
mechanism has been evidenced \cite{nagoshi}, but most importantly in 
laboratory experimental settings, tumour growth has been shown to be favoured by
disruptions of the normal circadian rhythm \cite{fili2005}.
{From} a mathematical point of view, cell population growth  is well described
by physiologically structured equations, see \cite{MD, chio}, and the
first eigenvalue of the underlying differential operator is the natural quantity that
accounts for the  growth of the system. Our purpose is to compare
the first eigenvalues
in the case of constant and periodic coefficients. It is simpler to 
consider in the first place
differential systems;  with periodic and nonnegative coefficients this
eigenvalue is nothing but  the (first) Floquet eigenvalue; with constant
coefficients it refers to the usual Perron eigenvalue.
Surprisingly, it is possible to prove in  great generality that
the Floquet eigenvalue is larger than the Perron eigenvalue with an
appropriate arithmetico-geometric average of the coefficients.
Precise statements, and proofs, are given in the first subsection in
the case of a differential system. The result also extends to discrete
time systems (section \ref{sec:dis}) and to Partial Differential
Equations and we give, in a fourth section, an application to
age-structured systems. Finally in a fifth section, we briefly 
comment on the relevance of these models to physiological systems.

\Sect{Differential systems}
Let $t\mapsto A(t)$ be a $T$-periodic 
map with values in $\R^{d\times d}$, integrable on $[0,T]$,
and let us consider
the differential equation $\dot X (t)= A(t)X(t)$,
with a prescribed initial condition $X(0)\in \R^d$
(when dealing with such differential equations, we will always
require $X$ to be absolutely continuous, and we understand that the
equality holds for almost all $t$).
We assume that for $i\neq j$, we have $A_{ij}(t)\geq 0$ for
almost every $t$, so that the flow in positive time of this
differential equation preserves the standard partial
ordering of $\R^d$. Hence, the spectral properties
of this differential equation belong to Perron-Frobenius
theory. In particular, the (first) Floquet
eigenvalue, $\lambdaper$, can be introduced by
means of the following positive Floquet problem:
there exists a $T$-periodic function $X$, with values
in $\R_+^d$ ($\R_+$ denotes the set of nonnegative real numbers),
non identically zero, such that
\be\label{eq:floquet}
\dot X(t) = A(t)X(t) -\lambdaper X(t)\enspace ,\;t\in \R   \enspace.
\ee
If $A$ satisfies some irreducibility properties, for instance
if $A_{ij}(t)>0$ for almost every $t$ and for $i\neq j$, $\lambdaper$
is uniquely defined by the previous property, and $X$ is unique
up to a multiplicative constant. However, the results of this note 
are also valid in the reducible case, in which the Floquet eigenvalue
can be defined as the maximal real number $\lambdaper$ such that there
exists a function $X$ with the above properties.

We now introduce the arithmetic mean of a $T$-periodic function $u(t)$ as
\be \<u>_a= \dis\frac 1 T \int_0^T u(s) ds
\label{e-ar}
\ee
and we define the constant
coefficient matrix  $\bar A$ with entries:
\be\label{eq:agaver}
\bar A_{ii} =\<A_{ii}>_a ,\;1\leq i \leq d,\qquad
\bar A_{ij}= \exp \left( \< \log(A_{ij})>_a \right) \qquad 
i\neq j,\; \; 1\leq i, j\leq d 
\enspace .
\ee
Observe that $\< \log(A_{ij})>_a\in \R\cup\{-\infty\}$
is well defined for $i\neq j$, because the positive
part of $\log A_{ij}(t)$ is integrable as soon as $A_{ij}(t)$
is integrable.
Since the off diagonal coefficients of the matrix $\bar A$ are nonnegative, 
we can apply the
Perron-Frobenius theory and  consider its first (Perron) eigenvalue,
$\lambda_s$, which is the maximal real number such that
there is a non-zero vector $U\in \R_+^d$
such that  $ \lambda_s U = \bar A U$. 
As discussed in Section~\ref{sec-relevance}, the following  result generalises
the one of \cite{clmp}.
\begin{theorem}
We have always 
$  \lambdaper  \geq  \lambda_s$.
\label{th:comp}
\end{theorem}
\proof We first prove the inequality when $A_{ij}(t)>0$, for
almost all $t$ and for $i\neq j$. Then, the function
$X$ in~\eqref{eq:floquet} is such that $X_i(t)>0$ for
all $t$ and for all $i$.
We set, 
\[
x_i(t):=\log X_i(t) \quad\text{ and  for } i\neq j,\quad \ell_{ij}(t):=\log A_{ij}(t) \enspace.
\]
{From} the differential system~\eqref{eq:floquet}, we obtain
\begin{align*}
\dot{x}_i(t)&=\sum_j X_i^{-1}(t)A_{ij}(t)X_j(t) -\lambdaper\\
&=\sum_{j\neq i}
\exp(-x_i(t)+\ell_{ij}(t)+x_j(t))+A_{ii}(t)-\lambdaper \enspace .
\end{align*}
Taking first the  arithmetic mean on  $[0,T]$ componentwise, and then
using Jensen's inequality, it comes for all $1\leq i \leq d$, 
\begin{align*}
0 &=\<\sum_{j\neq i} \exp(-x_i(t)+\ell_{ij}(t)+x_j(t))>_a
+\<A_{ii}(t)>_a-\lambdaper \enspace, \end{align*}
\begin{align*}
0 & \geq
\sum_{j\neq i} \exp(-\<x_i(t)>_a+\<\ell_{ij}(t)>_a+\<x_j(t)>_a)
+\<A_{ii}(t)>_a-\lambdaper \enspace. \end{align*}
Setting $\bar X_i:=\exp(\<x_i(t)>_a)$, with the definition~\eqref{eq:agaver} of
$\bar A$, this inequality also reads
\begin{align*}
0 & \geq
\sum_{j} \bar X_i^{-1}\bar A_{ij}\bar X_j-\lambdaper  \enspace,
\end{align*}
and thus,  multiplying by $\bar X_i$,
\[
 \bar A  \bar X \leq \lambdaper \bar X  \enspace . 
\]
Using the Collatz-Wielandt characterisation of the Perron eigenvalue
of $\bar A$, $\lambda_s= \min\{r; \; \exists Y\in \operatorname{int}\R_+^d, \bar A Y \leq r Y \}$,
we deduce that $\lambda_s\leq \lambdaper$. 
The general case is obtained by considering
the matrix with entries $A_{ij}(t)+\epsilon$,
with $\epsilon>0$, and by applying a continuity argument,
the details of which are left to the reader.
\qed

\Sect{Discrete systems: an inequality for the Perron eigenvalue of a geometric mean}
\label{sec:dis}

The same proof allows us to treat discrete systems. For $k\in \N$,
let $A(k)=A(k+p)$ be a $d\times d$, $p$-periodic matrix with nonnegative
coefficients. We define $\lambdaper$ to be the maximal real nonnegative
number such that there exists a non-zero $p$-periodic solution $X(k)$ 
with values in $\R_+^d$ to 
\be
\lambdaper \;  X(k+1) = A(k)X(k)  \enspace ,
\label{eq:discrete}
\ee
so that $\lambdaper^p$ is the Perron eigenvalue of the product $A(p-1)\cdots A(0)$.
We now set $\<u>_a:=p^{-1}(u(0)+\cdots +u(p-1))$ for
all $p$-periodic functions $u$, and we define the constant coefficient matrix
$\bar A$, with entries:
\be\label{eq:agaverds}
\bar A_{ij}=
\exp \left( \< \log(A_{ij})>_a \right) ,\; \;1\leq i,j\leq d
\ee
(unlike in the continuous time case, the nature
of the mean is the same for diagonal and off diagonal entries).
We denote by $\lambda_s$ the Perron eigenvalue
of the nonnegative matrix $\bar A$.
\begin{theorem}
We have again $\lambdaper \geq  \lambda_s$.
\label{th:compds}
\end{theorem}

\proof Since the Perron eigenvalue is a continuous
function of the entries of a matrix,
it suffices to show the inequality when the entries
$A_{ij}(k)$ are all positive. Then, the vectors $X(k)$
above have positive entries.
We set
\[
x_i(k):=\log X_i(k),\quad
\ell_{ij}(k):=\log A_{ij}(k), \quad \mu := \log(\lambdaper)
\enspace .  
\]
We take componentwise the logarithm in 
$\lambdaper X(k+1)= A(k)X(k)$ and arrive
at
\[
x_i(k+1)-x_i(k)=
\log\left(\sum_{j}
\exp\big(-x_i(k)+\ell_{ij}(k)+x_j(k)\big)\right)-\mu \enspace .
\]
Taking the  arithmetic mean for $k=0,\ldots,p-1$, it comes
\[
0 = \<\log\left(\sum_{j}
\exp\big(-x_i(k)+\ell_{ij}(k)+x_j(k)\big)\right)>_a -\mu \enspace. 
\]
Next, because the function $f(y_1,\ldots,y_d)=\log(\sum_j \exp(y_j))$
is convex, we apply Jensen's inequality and obtain 
\[
0  \geq  \log\left(\sum_{j}
\exp\big(-\<x_i(k)>_a+\<\ell_{ij}(k)>_a+\<x_j(k)>_a\big)\right)
-\mu \enspace, 
\]
and, exponentiating, we get
\[
\lambdaper \bar X \geq \bar A \bar X \enspace,
\]
with $\bar X_i:=\exp(\<x_i(k)>_a)$.
Using again the Collatz-Wielandt characterisation of the Perron eigenvalue
of $\bar A$, we deduce that 
$\lambdaper\geq \lambda_s$.
\qed

\Sect{An age-structured system for the cell division cycle}

General references and experimental
validations on the topic of structured population dynamics and
cell cycle can be found in \cite{chio,MD}. For a  recent mathematical
approach based on entropy properties, we refer to \cite{michel1,
perthamebook}. Here and following earlier work \cite{clmp}, we model
our population of cells by a Partial Differential Equation for the
density $n_i(t,x)\geq 0$ of cells with age $x$ in the phase
$i=1,\ldots,$ at time $t$,
\begin{equation} \left\{ \begin{array}{l} 
\frac{\partial}{\partial t} n_i(t,x)+ \frac{\partial}{\partial
x}n_i(t,x)+[d_i(t,x) + K_{i\to i+1}(t,x)] n_i(t,x) =0  \enspace ,
\\ [2mm]
n_i(t,x=0)= \dis \int_{x'\geq 0}K_{i-1\to i}(t,x')\; n_{i-1} (t,x')\;
dx',
\quad 2\leq i \leq I,
\\ [3mm]n_1(t,x=0)=2 \dis \int_{x'\geq 0}K_{I \to 1}(t,x')\; n_{I}
(t,x')\; dx'  \enspace .\end{array} \right. \label{eq:cellcycle}
\end{equation}
Here and below we identify ${I+1}$ to $1$. We have denoted by
$d_i(t,x)\geq 0$ the apoptosis rate, by $K_{i\to i+1}$ the transition
rates from one phase to the next, and the last one ($i=I$) is mitosis
where the two cells separate. These coefficients can be constant in
time (no circadian control) or time $T$-periodic in order to take
into account the circadian rhythm. Our assumptions are 
\begin{equation}
K_{i\to i+1}(t,x) \geq 0, \; d_i(t,x) \geq 0 \qquad \text{are
bounded},\label{eq:as1} 
\end{equation}
and, setting \[
\dis\min_{0\leq t \leq T} K_{i\to i+1}(t,x)
:= k_{i\to i+1}(x), \quad \dis\max_{0\leq t \leq T} [d_i+ K_{i\to
i+1}]:=\mu_i(x), \quad M_i(x)= \dis\int_0^x\! \mu_i(y) dy,
\] 
\begin{equation} 
\dis\prod_{i=1}^{I}\dis\int_0^\infty k_{i\to i+1}(y)e^{-M_i(y)} dy
>1/2\enspace .
\label{eq:as2} 
\end{equation}
With these assumptions and following \cite{michel1}, one can again
introduce the growth rate (Floquet eigenvalue)  of the system:
$\lambdaper \in \R$ such that there is a unique $T$-periodic {\em
positive} solution to the system:
\begin{equation} 
\left\{ \begin{array}{l}
 \frac{\partial}{\partial t}N_i(t,x)+ \frac{\partial}{\partial
x}N_i(t,x)+[d_i(t,x) +\lambdaper +K_{i\to i+1}(t,x)] N_i(t,x) =0,
 \\ [2mm] 
N_i(t,x=0)= \dis\int_{x'\geq 0}K_{i-1\to i}(t,x')\; N_{i-1} (t,x')\;
dx'\enspace ,
\quad 2\leq i \leq I \enspace ,
\\ [2mm]
N_1(t,x=0)=2 \dis\int_{x'\geq 0}K_{I \to 1}(t,x')\; N_{I} (t,x')\;
dx' \qquad \dis\sum_{i=1}^I\dis\int_{x\geq 0} N_i(t,x) dx=1\enspace  .
\end{array} \right.
\label{eq:per}
\end{equation}

As in formula \fer{eq:agaver}, we can define the averages
\be\label{eq:avercc}
\left\{ \begin{array}{l}
\<d_i(x)>_a =\dis \frac 1 T \dis \int_0^T d_i(t,x)dt, \qquad  \<K_{i\to
i+1}(t,x)>_a  = \frac 1 T\dis  \int_0^T  K_{i\to i+1}(t,x) dt\enspace
,
\\ [2mm]
\<K_{i\to i+1}(t,x)>_g = \exp\left( \dis \frac 1 T \int_0^T  \log
\big(K_{i\to i+1}(t,x) \big) dt \right)\enspace . 
\end{array} \right.
\ee
These averages define the Perron eigenvalue $\lambda_s \in \R$, which 
is such that there is a  unique positive solution to the system 
\be
\left\{ \begin{array}{l}\
\frac{\partial}{\partial x} \bar{N}_i(x) +[\<d_i(x)>_a+\lambda_s+
\<K_{i\to i+1}(t,x)>_a]\bar N_i=0\enspace ,
\\ [2mm]
\bar N_i(x=0)=\dis \int_{x'\geq 0} \<K_{i-1\to i}(t,x')>_g
\bar{N}_{i-1}(x')d x',\; i\neq 1\enspace ,
\\ [3mm]
\bar N_1(x=0)=2 \dis \int_{x'\geq 0} \<K_{I\to 1}(t,x')>_g
\bar{N}_{I}(x')d x' \enspace.
\end{array} \right.
\label{e-g3} \ee
We have the following analogue of Theorem \ref{th:comp}.
\begin{theorem}
Under assumptions \fer{eq:as1}--\fer{eq:as2}, we still have
$\lambdaper  \geq  \lambda_s$.
\label{th:compcc}
\end{theorem}
Observe that the arithmetic mean of the coefficients $K_{i\to i+1}$
is taken in the PDEs, whereas their geometric mean is taken in the
integral equations. Hence, an artificial loss rate of cells
$\<K_{i\to i+1}>_a-\<K_{i\to i+1}>_g$ from phase $i$ to phase $i+1$
arises in the averaged model.

\proof The proof differs slightly from that of Theorem \ref{th:comp}
because, working with  $x \in \R_+$ the corresponding quantities
cannot always be normalized as measures. Therefore, following
\cite{clmp}, we define $q_i(x)= \< \log N_{i}(t,x)-\log \bar N_{i}(x)
>_a$ and we have, up to the insertion of a factor $2$ when $i=1$,
\[
q_i(x=0)=\<\log(\frac{N_i(t,0)}{\bar{N}_i(0)})>_a
= \<\log\int \<K_{i-1\to i}(x)>_g
\frac{\bar{N}_{i-1}(x)}{\bar{N}_i(0)}
\frac{N_{i-1}(t,x)}{\bar{N}_{i-1}(x)} \frac{K_{i-1\to
i}(t,x)}{\<K_{i-1\to i}(x)>_g}dx>_a\enspace . 
\]
We can now define the probability measures
\[
d\mu_i(x)=\<K_{i-1\to i}(x)>_g \frac{\bar{N}_{i-1}(x)}{\bar{N}_i(0)}
\quad 
\text{ for }i\neq 1,
\qquad 
d\mu_1(x)=2\<K_{I\to 1}(x)>_g
\frac{\bar{N}_{I}(x)}{\bar{N}_1(0)}  \enspace .
\]
Using Jensen's inequality, we obtain 
\begin{align*}
q_i(x=0)&\geq \< \int \log(\frac{N_{i-1}(t,x)}{\bar{N}_{i-1}(x)}
\frac{K_{i-1\to i}(t,x)}{\<K_{i-1\to i}(x)>_g}) d\mu_i(x)>_a
\\ 
&= \< \int \log(\frac{N_{i-1}(t,x)}{\bar{N}_{i-1}(x)})d\mu_i(x)>_a
+\< \int \log(\frac{K_{i-1\to i}(t,x)}{\<K_{i-1\to i}(x)>_g})
d\mu_i(x)>_a
\\ &=
\int q_{i-1}(x)d\mu_i(x) \enspace,
\end{align*}
since by definition of $\<K_{i-1\to i}(x)>_g$ we have $\dis\int \<
\log(\frac{K_{i-1\to i}(t,x)}{\<K_{i-1\to i}(x)>_g})>_a d\mu_i(x)=0$.
Because $q_i$ satisfies
\[
\frac{\partial}{\partial x} q_i + \lambdaper -\lambda_{s} =0 \enspace,
\]
we arrive at
\[
q_i(x=0)  \geq \dis \int q_{i-1}(x) \; d\mu_i(x) = \dis \int
[q_{i-1}(0)+(\lambda_{s} - \lambdaper) x] \; d\mu_i(x)\enspace .
\]
Therefore, summing over $i$ from $1$ to $I$, we have obtained the
result since
$$
0 \geq (\lambda_{s} - \lambdaper) \dis \sum_{i=1}^I \int x \;
d\mu_i(x) \enspace .
$$
\qed
\Sect{Physiological relevance of the model to the question of 
circadian control on tumour growth}\label{sec-relevance}
In \cite{clmp}, with different averaging of transition rates
$K_{i\to i+1}$, we had found no natural hierarchy between 
$\lambdaper$ and $\lambda_s$, whereas periodic control exerted on the 
sole apoptosis rates $d_i$ produced a Floquet eigenvalue which was 
shown to be always higher than the corresponding Perron eigenvalue obtained 
by arithmetic averaging of the $d_i$.  The present note gives natural 
hypotheses that are needed to compare the eigenvalues and conclude
that $\lambdaper>\lambda_s$ when the coefficients $K_{i\to i+1}$
also vary. {From}
these results a question arises: if $\lambdaper>\lambda_s$, then 
how come that experimental results \cite{fili2005} show faster tumour 
growth when the normal circadian rhythm is disrupted by irregular 
light inputs? Experimental tumour growth curves with and without 
jet-lag-like circadian clock disruption \cite{fili2005}
essentially differ at the beginning of proliferation, when
it shows almost pure exponential behaviour, so that 
possible discrepancies between theoretical and experimental results are most 
likely not due to not taking into account nonlinear feedback. Clearly, 
setting at a constant value (instead of periodic) the coefficients $d_i$ or $K_{i\to i+1}$ 
is not sufficient to take into account the complexity involved in the disruption 
of circadian control. To this purpose, it may be necessary to represent 
more complex control mechanisms involving, e. g., the inhibition of cyclin dependent
kinases by clock-controlled genes such as Wee1 \cite{gol91} together with 
elaborate models of the molecular circadian clock \cite{lel04} and its disruptions.

\end{document}